\documentclass[11pt]{amsart}

\usepackage{amssymb}

\usepackage{amsmath,latexsym,amsthm}

\usepackage{amsfonts}
\usepackage[greek,english]{babel}

\newtheorem{claim}{Claim}
\newtheorem*{theorem*}{Theorem}
\newtheorem*{corollary*}{Corollary}

\def\bel{\begin{equation}\label}
\def\eeq{\end{equation}}

\newtheorem{remark}{Remark}[section]
\newtheorem{definition}{Definition}[section]

\newtheorem{theorem}{Theorem}[section]

\newtheorem{proposition}[theorem]{Proposition}

\newtheorem{example}[theorem]{Example}

\newcommand\epf{\hfill{$\square$}\medskip}

\def\ds{\displaystyle}

\def\bega{\begin{array}}
\def\enda{\end{array}}

\def\bepmatrix{\begin{pmatrix}}
\def\enpmatrix{\end{pmatrix}}

\def\bel{\begin{equation}\label}
\def\eeq{\end{equation}}
\newcommand\ee{\end{equation}}
\def\benl{\begin{equation*}}
\def\eenl{\end{equation*}}

\def\be{\begin{equation}}
\def\beq{\begin{equation}}
\def\bel{\begin{equation}\label}
\def\eeq{\end{equation}}

\newcommand\ba{\begin{array}}
\newcommand\ea{\end{array}}

\def\begi{\begin{itemize}}
\def\endi{\end{itemize}}

\newcommand{\cR}{\mathbb{R}}
\newcommand{\cQ}{\mathbb{Q}}
\newcommand{\cN}{\mathbb{N}}

\def\C{\mathcal{C}}

\def\I{\mathcal{I}}

\def\L{\mathcal{L}}

\def\T{\mathcal{T}}
\def\U{\mathcal{U}}

\def\ub{\bar{u}}

\def\xb{\bar{x}}

\def\eps{\varepsilon}

\def\ubf{{\varphi}}
\def\vbf{{\psi}}
\def\ybf{{ y}}

\title[Limit solutions]{$\L^1$ limit solutions for control systems$^1$}

\begin{document}

\author[M.S. Aronna]{M. Soledad Aronna}
\address{M.S. Aronna, IMPA, Estrada Dona Castorina 110, Rio de Janeiro 22460-320, Brazil}
\email{aronna@impa.br}

\author[F. Rampazzo]{Franco Rampazzo}
\address{F. Rampazzo, Dipartimento di Matematica ,
Universit\`a di Padova\\ Padova  35121, Italy}
\email{rampazzo@math.unipd.it}

\maketitle

\footnotetext[1]{Article published in  J. Differential Equations, 258(3): 954-979, 2015}

\begin{abstract}
For a control Cauchy problem
$$\dot x= {f}(t,x,u,v) +\sum_{\alpha=1}^m g_\alpha(x) \dot u_\alpha,\quad x(a)=\bar x,  $$
on an interval $[a,b]$,
we propose  a notion of  {\em limit solution} $x,$ verifying  the following properties: i) $x$ is defined  for $\L^1$ {(impulsive)}   inputs  $u$  and for standard, bounded measurable, controls  $v$; ii) in the commutative case (i.e. when $[g_{\alpha},g_{\beta}]\equiv 0,$ for all $\alpha,\beta=1,\dots,m$), $x$ coincides with the solution one can obtain via the  change of coordinates that makes the $g_\alpha$ simultaneously  constant;
iii) $x$ subsumes  former concepts of solution valid for the generic, noncommutative case.
 In particular,  when $u$ has  bounded variation, we investigate the relation between  limit solutions and (single-valued)  graph completion solutions.
Furthermore, we prove consistency with the classical Carath\'eodory solution when $u$ and $x$ are absolutely continuous.

Even though some specific problems are better addressed by means of special representations of the solutions,  we believe that various theoretical issues call  for a unified notion of trajectory.  For instance, this is the case of optimal control problems, possibly with state and endpoint constraints, for which no extra assumptions (like e.g. coercivity, bounded variation, commutativity) are made in advance.

\vspace{6pt}


\vspace{6pt}

\end{abstract}


\section{Introduction}

\subsection{Limit solutions}

Because of the presence of the derivative $\dot u$,  control systems of the form \begin{align}
 \label{E}\tag{E} \dot x &= {f}(t,x,u,v) +\sum_{\alpha=1}^m g_\alpha(x) \dot u_\alpha,\quad t\in [a,b], \\
\label{IC}\tag{IC} x(a)&=\xb\end{align} are sometimes called {\it impulsive}. This is connected with the  need  of implementing non regular, possibly discontinuous, controls $u$ --a need raised e.g. by  lack of coercivity (in $\dot u$) in a minimum problem. Actually,
the issue of giving a notion of solution for a control  system like (E) becomes non standard as soon as $u$ is not absolutely continuous (while $v$ can  be an ordinary bounded measurable control). As is well-known, the core of the question  resides in  the interaction between the $x$-dependence of the vector fields $g_1,\dots, g_m$ and the unboundedness of the derivative $\dot u.$ 
Let us point out that several physical control settings (\cite{BreAldo89,Ram90,Mar91,Bra91,Cat08,GajRamRap08}) are naturally modeled by an equation like \eqref{E} --possibly, as in some mechanical applications, \eqref{E} being a first order reduction of a higher order equation.

Notions of output are well-established at least in two situations:
\begin{itemize}
\item[(i)] The case of  {\it commutative systems}, characterized by the triviality of the associated Lie algebra. This means that  $[g_\alpha,g_\beta]\equiv 0,$ for all $\alpha,\beta=1,\dots,m$ (see e.g. \cite{BreRam91,Sar91,Dyk94,AroRam14}).  Actually, under such hypothesis, a  notion of solution to \eqref{E}-\eqref{IC} has been established  for controls  $u\in \L^1$ (here we use $\L^1$ to denote the set of {\em pointwise defined} Lebesgue integrable functions, while $L^1$  refers to the usual quotient with respect to the Lebesgue measure), possibly with {\it unbounded variation}. These solutions, whose definition essentially relies on a state space diffeomorphism  induced by the family of vector fields $\{g_1,\dots, g_m\},$   are {\it pointwise defined} and verify nice properties of well-posedness.
\item[(ii)] The case when {\it  the controls $u$ have   bounded variation} (and the Lie algebra is allowed to be non trivial).     In order to  verify standard  robustness properties, this  notion of solution requires a specification on how the discontinuities of $u$  are bridged.  This type of solutions (see e.g.  \cite{Ris65,BreRam88,DalRam91,SilVin96,PerSil00,CamFal99,Kar06,WolZab07,AruKarPer11}), which are not single-valued at the jump instants\footnote{But a single-valued version is  proposed here, via the notion of {\em clock} (see Definition \ref{DefClock}).},  are  described by different authors in fairly equivalent ways,  and will be called here {\it graph completion solutions.} 
\end{itemize}

To our knowledge, no available definition of solution is applicable to  the union of the above-described subclasses of equations.
 On the other hand,
let us  point out that  a general  notion of solution (including  the above cases but not limited to them) is  of interest for those questions  where specific properties (of either  the controls or the vector fields)   are not  known {\it a priori}.

The main objective of the present paper consists in the discussion of a unified  concept  of solution to \eqref{E}-\eqref{IC}.
As a matter of  fact, we will  adopt a {\it density} characterization valid in the commutative case as  a definition of solution  for the general case. This concept will be denominated {\it ``limit solution"}.

Let us give immediately the  notion of limit solution, postponing the precise specification of  the regularity and growth assumptions on the vector fields.

\begin{definition}
\label{edsdef}
Consider an initial value $\xb \in \cR^n,$ controls $u \in \L^1([a,b];U)$  and $v\in L^1([a,b];V).$
\begin{enumerate}
\item {\sc (Limit Solution)} We say that a map $x:[a,b]\to \cR^n$ is a {\em limit solution} of the Cauchy problem \eqref{E}-\eqref{IC} corresponding to $\xb$ and $(u,v)$ if, {\em for every $\tau \in [a,b],$}  there exists a sequence of absolutely continuous controls $(u^\tau_k)$ from $[a,b]$ into $U$  such that:
\begin{itemize}
\item[{(i$_\tau$)}] for each $k \in \cN,$ there exists a (Carath\'eodory) solution
$$
x^\tau_k:[a,b]\to\cR^n
$$
of  \eqref{E}-\eqref{IC} corresponding to the initial value $\xb$ and the control $(u^\tau_k, v),$ and the sequence $(x^\tau_k)$  is bounded;
\item[{(ii$_\tau$)}] one has
$$
|(x^\tau_k,u^\tau_k)(\tau)- (x,u)(\tau)|+ \|(x^\tau_k,u^\tau_k)- (x,u)\|_1\to 0,
$$
where $\|\cdot\|_1$ denotes the $L^1$-norm.
\end{itemize}
\item {\sc (Simple limit solution)}
We say that  a limit solution $x:[a,b]\to \cR^n$ is  {\em simple} if $(u^\tau_k)$ can be chosen independently of $\tau.$ In this case we write $(u_k)$ to refer to the approximation sequence.
\item {\sc (BV simple limit solution)}\label{sedsdef}
We say that a simple  limit solution $x:[a,b]\to \cR^n$ is a {\em BV (bounded variation) simple limit  solution} if the approximating inputs  $u_k$ have equibounded variation.
\end{enumerate}
\end{definition}
The idea of defining control-trajectories pairs as limits of regular trajectories is obviously not new (we refer to e.g. \cite{BreRam91,LiuSus91,Sar91} and the  references therein). 
Let us insist on the fact that we are interested in integrable, possibly discontinuous inputs $u$ \footnote{If attention were confined to  the case of {\it continuous} inputs $u,$ the subtlest instrument of investigation might be T. Lyons' {\it rough paths} (\cite{Lyo94,Lyo02}), which   are introduced --via iterated integrals-- as metric extensions of continuous input-output maps  (see  paragraph \ref{continuousu}).}.
 Furthermore, the fact that in the definition of limit solution the  approximating  sequences $(u_k^\tau)$ do depend on the considered  time  $\tau$    is essential if we aim to a notion of everywhere defined  solution corresponding to  $\L^1$ controls $u$ (see Example \ref{ex-suss} below). Notice incidentally that such dependence makes the topology underlying the definition of limit solutions  not metrizable.

In the presence of $u$-jumps and unless the Lie algebra generated by the $g_\alpha$'s is trivial,  an irremediable  multiplicity of solutions has to be considered. Nevertheless, the following facts hold true:
\begin{itemize}
 \item[a)]{\it (Generality)} {The notion of limit solution  subsumes the two  concepts of solutions described above and seems to fit natural generalizations such as inputs with unbounded variation, including those which are not pointwise limit of continuous maps (see Example \ref{ex-suss} below).}
 \item[b)] {\it (Representability of BV simple limit solutions as graph completion solutions)} When $u$ has bounded variation, the notions of graph completion solution  and that of BV simple  limit solution turn out to be equivalent (Theorem \ref{gcpdgen1}).
 \item[c)]  {\it (Consistency with the classical case)} As we have already remarked, if $u$ is absolutely continuous, any absolutely continuous BV simple  limit solution $x$  corresponding to $u$ and to a given $v$  coincides with the classical Carath\'eodory solution (Theorem \ref{cons}). However it is worth noticing that smooth limit solutions may exist without being Carath\'eodory solutions, as shown in Example \ref{ex-suss} below.
     \end{itemize}

Here is a brief outline of the paper.
 Section \ref{limit-def} is opened by an example that aims to provide some justification for the introduction of such general solutions. Then we   discuss properties of  uniqueness and consistency (with  Carath\'eodory solutions). Successively, existence and uniqueness are  proved in the commutative case. Finally, the existence of BV simple  limit solutions associated to controls with bounded variation is stated (the proof being provided in Section \ref{proofSec}).
In Section \ref{commSec} we address the  notion of limit solution in the case where the vector fields $g_1,\dots,g_m$ commute (see \cite{BreRam91,AroRam13,AroRam14}, and also  \cite{Dyk94,Sar91}).
We observe that, unlike the noncommutative case,  uniqueness and  continuous dependence is obtained also for discontinuous, $\L^1$ controls $u.$ 
Section \ref{noncommSec} deals with the noncommutative case.  After introducing (a single-valued version of) the notion of graph completion solution, we prove a result which shows how this concept can  be embedded in the general notion of limit solution. Precisely: a trajectory $x$ is  a graph completion solution if and only if it  coincides {pointwise} with a BV simple  limit solution (Theorem \ref{gcpdgen1}). By proving that for each $u$ with bounded variation a graph completion solution actually exists, one gets  the existence result of BV simple  limit solutions stated in Section \ref{limit-def}. This proof is placed in Section  \ref{proofSec}. Incidentally, as a byproduct, one obtains  that every function with bounded variation (and values on a set with the Whitney property) can be pointwise approximated by means of (absolutely) continuous functions with {\em equibounded variation} (see Theorem \ref{whi}).
In Section \ref{conclSec} we briefly illustrate some open issues concerning, in particular, the case where the vector fields $g_\alpha$ depend on $v$ as well, and the case where these fields are not differentiable.


\subsection{Some notation}

Let $I,$ $E$ be a closed interval and a subset of the Euclidean  space $\cR^d,$ respectively. For any $L\in [0,+\infty[,$ we let $ BV_L(I;E)$ be the set of functions $h:I\to E$ with total variation bounded by $L,$ and we set $BV(I; E)=\bigcup_{L\geq0} BV_L(I;E).$
We write ${\rm Var}_I(h)$ to refer to the variation of $h$ in the interval $I.$
Moreover, we shall recall that $\L^1([a,b];E)$ denotes the set of pointwise defined Lebesgue integrable functions from $I$ to $\cR^d$ with values in $E,$ while with $L^1([a,b];E)$  we indicate the corresponding  family of equivalence classes (with respect to the Lebesgue measure). With $AC(I;E)$ we refer to the set of absolutely continuous maps from $I$ to $E,$ and we  let $AC_L(I;E)\subset AC(I;E)$ represent the subset made of the  absolutely continuous  functions having variation bounded by $L,$ namely $AC_L(I;E):=  AC(I;E)\cap  BV_L(I;E).$
The notation ${\rm Lip}(I;E)$ will indicate the set of Lipschitz continuous functions from $I$ to $E.$  $\C^k(\cR^n;\cR^d)$ will denote the space of $k-$times continuously differentiable functions defined on $\cR^n,$ with values in $\cR^d.$ The elements in $\C^k(\cR^n;\cR^d)$ will be sometimes  called {\em functions of class $\C^k.$}

We will say  that a real function $\varphi_0:[a,b]\to \cR$ is {\em increasing} (respectively, {\em decreasing}) if for every pair $t_1,t_2 \in [a,b]$ verifying $t_1 < t_2,$ one has $\varphi_0(t_1)\leq \varphi_0(t_2)$ (respectively, $\varphi_0(t_1)\geq \varphi_0(t_2)$). We use {\em strictly increasing} (respectively, {\em strictly decreasing}) when the corresponding inequality is strict.


\subsection{Structural hypotheses}
Throughout the paper we shall assume the following hypotheses on the control sets $U,$ $V$ and the functions $f,g_1,\dots,g_m:$
\begin{itemize}
\item[(i)]
  $U$ is a compact subset of $\cR^m$ such that, for every bounded interval $I \subset \cR,$ for each $\tau \in I,$ and for every  function $u\in\L^1(I;U),$ there exists a sequence $(u^\tau_k) \subset AC(I;U)$ such that
    $$
 |u_k^\tau(\tau)- u(\tau)| + \|u^\tau_k-u\|_1 \to 0,\quad \text{when $k \to \infty.$}
  $$
  (Convex sets verify this hypothesis. Furthermore, it is not difficult to show that also the closure of any open connected, bounded subset with a Lipschitz boundary meets (i)).
  \item[(ii)]
The set $V \subset \cR^l$ is compact.
 \item[(iii)] For each $(x,u,v)\in \cR^n\times \cR^m \times V,$ the map $t\mapsto f(t,x,u,v)$ is measurable on $[a,b];$
 for each $t\in [a,b],$ the map $(x,u,v) \to f(t,x,u,v)$ is continuous on $\cR^n \times \cR^m \times V$ and, moreover,  the map
 $$
 (x,u) \mapsto f(t,x,u,v)
 $$
 is locally Lipschitz on $\cR^n \times \cR^m,$ uniformly in $(t,v) \in [a,b] \times V.$
 \item[(iv)] The vector fields $g_\alpha:\cR^n \to \cR^n$ are of class $\C^1.$
  \item[(v)] There exists $M>0$ such that
  $$
  \left| \Big(f(t,x,u,v),g_1(x),\dots,g_m(x)\Big)\right| \leq  M(1+|(x,u)|),$$ for every $(x,u) \in \cR^n\times \cR^m$ uniformly in $(t,v) \in [a,b]\times V.$\end{itemize}

Notice that, under conditions (ii)-(v) above, for any initial value $\xb \in \cR^n$ and each control pair $(u,v)\in AC([a,b];U)\times {L^1}([a,b];V),$ there exists a unique Carath\'eodory solution to the Cauchy problem \eqref{E}-\eqref{IC}. We let $x[\xb,u,v]$ denote this solution.

\subsubsection{Commutativity and completeness}
For some results we shall assume the hypothesis described below.

\if{
Set $\tilde g_\alpha:= \sum_{j=1}^n g_{\alpha,j}\frac{\partial}{\partial x_j} + \frac{\partial}{\partial u_\alpha},$ for every $\alpha=1,\dots,m,$ where
$$
\left( \frac{\partial}{\partial x_1},\dots,\frac{\partial}{\partial x_n},\frac{\partial}{\partial u_1},\dots,\frac{\partial}{\partial u_m}\right)
$$
is the canonical basis of $\cR^{n+m}.$
We use $[\tilde {g}_\alpha,\tilde {g}_\beta] $  to denote  the {\em Lie bracket}  of the vector fields $\tilde {g}_\alpha$ and $\tilde {g}_\beta.$  In coordinates, it is  defined by
$$
\left[\tilde {g}_\alpha,\tilde {g}_\beta \right] =
D\tilde{g}_\beta \, \tilde g_\alpha - D\tilde {g}_\alpha \, \tilde g_\beta.
$$
Here $D$ denotes  the differential operator in the  $(x,u)$ state-space, so that, in particular, the last $m$ components of $\left[\tilde {g}_\alpha,\tilde {g}_\beta \right] $ are zero.
}\fi

\begin{definition}
\begin{itemize}
\item[(i)] We say that $g_1,\dots,g_m$ verifies the {\em commutativity hypothesis}  if
$$
\left[{g}_\alpha,{g}_\beta \right](x)=0,\quad \forall x \in \cR^{n},
$$
where $[g_\alpha,g_\beta]:= D{g}_\beta \,  g_\alpha - D {g}_\alpha \,  g_\beta$ is the {\em Lie bracket} of $g_\alpha$ and $g_\beta.$
\item[(ii)]
The vector field $ g_\alpha$ is {\em complete} if the solution of the Cauchy problem
$$
\dot x(t)=g_\alpha(x(t)),\quad x(0)=\xb,
$$
is (uniquely) defined on  $\cR.$
\end{itemize}
\end{definition}

The following assumption will be adopted for some results of this article.

\vspace{5pt}

\noindent {\bf Hypothesis (CC):} We say that the functions $g_1,\dots,g_m$ satisfy the {\em Hypothesis (CC)} if they verify the commutativity hypothesis and the vector fields $ g_1,\dots, g_m$ are complete.

\section{Limit solutions}
\label{limit-def}

\subsection{An example}

Let us consider the Cauchy problem \eqref{E}-\eqref{IC}.
In the Definition \ref{edsdef} of limit solution, the fact that the choice of $(u^\tau_k)$ depends on the time $\tau$  might appear awkward. However, this feature is essential to incorporate the case $u\in\L^1([a,b];U),$ as shown in Example \ref{ex-suss} below.  On the other hand, some important results are valid for the case where the choice of  $(u_k^\tau)$ is actually independent of the point $\tau,$ which explains why we have introduced the notion of (possibly BV)  simple  limit solution.

\begin{example}\label{ex-suss}
{\rm
 Let us set  $h(x):=(x_3-x_5)^2 + x_1^2 +x_2^2 $ and define the vector fields   on $\cR^6$
$$
 f(x):=-\frac{\partial}{\partial x_5}+\frac{\partial}{\partial x_6}h(x),
$$
\be
\label{exg}
g_1:=\frac{\partial}{\partial x_1}+x_2\frac{\partial}{\partial x_3},\quad 
g_2:=\frac{\partial}{\partial x_2}-x_1\frac{\partial}{\partial x_3},\quad g_3:=x_4\frac{\partial}{\partial x_4}.
\ee
 Given a map $t\mapsto\phi(t)$ belonging to $\L^1([0,1];\cR)$ with $\phi(0)=0,$
 and a finite  subset $I\subset ]0,1]$, let us  consider the optimal control problem 
\be
\tag{$P_{\phi}$}
\inf \, \mathcal I_\phi(x),
\ee 
where 
$$
\mathcal I_\phi(x):=x_6(1) + \sup_{t\in I} \left(x_4(t)-e^{\phi(t)}\right)^2,
$$ 
$x(\cdot)$ being a solution of 
\benl
\begin{split}
& \dot x = f(x) + g_1(x)\dot{u}_1 +g_2(x)\dot{u}_2 + g_3(x)\dot{u}_3,\\
& x(0)=(0,0,1,1,1,0),
\end{split} 
\eenl
for some input $u(\cdot).$ 
                         
Let us observe that the Lie bracket $[g_1,g_2]$ is not  equal to zero. This rules out the possibility of finding a trivializing coordinates' change which would allow for $\L^1$ controls $u$ and uniqueness of the corresponding solutions (according to  \cite{BreRam91,AroRam14}). 

To state Facts 1 and 2 below, let us 
say that {\bf Hypothesis (S)} holds for a function $\phi:[0,1]\to\cR$ if there exists a sequence $(\phi_k)$ of absolutely continuous equibounded maps converging to $\phi$ pointwise everywhere in $ [0,1],$ and such that $\|\phi_k-\phi\|_1 \to 0.$

\noindent {\bf Fact 1.} {\em Assume that Hypothesis (S) holds true for some function $\phi,$ then the Carath\'eodory solutions $x_k$ corresponding to the (absolutely continuous) controls
$$
u_k(t): = \left( \frac{\cos{kt}-1}{\sqrt{k}},\frac{\sin{kt}}{\sqrt{k}},\phi_k(t) \right),
$$
verify
\be
\label{limitxk}
 x_k(t)\to\tilde x(t):=(0,0,1-t,e^{\phi(t)},1-t,0),\qquad\forall t \in [0,1].
\ee
}
Indeed, each $x_k$ is given by
\benl
x_k(t) =  \left(u_{k,1}(t),u_{k,2}(t),1-t+\frac{\sin{kt}}{k},\frac{e^{\phi_k(t)}}{e^{\phi_k(0)}},1-t, a_k(t) \right),
\eenl
with $a_k(t):=\left( \frac{1}{2k^2} + \frac{2}{k} \right)t -\frac{\sin 2kt}{4k^3}-\frac{2\sin kt}{k^2}$,
which implies \eqref{limitxk}. Moreover, setting $\tilde u(t):=(0,0,\phi(t)),$ one gets that
$\big( (u_k,x_k) \big)$ converges to $ \big( \tilde u(t)
, \tilde{x}(t) \big)$  pointwise everywhere and in the $L^1$-topology.
In particular, one has 
\bel{zero}
{\mathcal I}_\phi(x_k)\to 0 =
{\mathcal I}_\phi(\tilde x).
\eeq
Hence,  $\tilde x(\cdot)$ is a simple limit solution and
\bel{zero1}
\inf \, {\mathcal I}_\phi(x) =
{\mathcal I}_\phi(\tilde x),
\eeq
where the infimum on the left-hand side is taken over Carath\'eodory solutions (corresponding to absolutely controls $u$).

\vspace{5pt}
 
\noindent {\bf Fact 2.} {\em Let us now discard Hypothesis (S) on $\phi.$  Consider, for   every $\tau\in[0,1],$ a sequence $(\phi_k^\tau)$ --surely existing-- of absolutely continuous equibounded maps converging to $\phi$  in the $L^1$-norm and verifying $\phi_k^\tau(\tau) \to \phi(\tau).$  Then the Carath\'eodory solutions $x_k^\tau$ corresponding to the absolutely continuous controls
$$
u_k^\tau(t): = \left( \frac{\cos{kt}-1}{\sqrt{k}},\frac{\sin{kt}}{\sqrt{k}},\phi_k^\tau(t) \right),
$$
verify 
$$
\left | (u_k^\tau,x_k^\tau)(\tau)-\big( \tilde u,\tilde x \big)(\tau)
\right| + \left\| (u_k^\tau,x_k^\tau)(\cdot)-\big( \tilde u,\tilde x \big)(\cdot)
\right\|_1\to 0.
$$}
Hence, even though Hypothesis (S) is not standing,  $\tilde x(\cdot)$ is a (not simple) limit solution associated to $\tilde u(\cdot).$
Note that, however, $\tilde x(\cdot)$ is not a Carath\'eodory solution, even if $\phi$ were assumed to be absolutely continuous.

\vspace{8pt}

Some considerations are at issue. 
First,  Fact 1 can be viewed as a motivation for the introduction of the notion of {\em simple limit solution,} possibly corresponding to sequences of inputs with variations tending to infinity.
Secondly, Fact 2  shows why in order to represent weak solutions to minimum problems it may prove important to enlarge  the class of simple limit solutions to 
 maps which are not pointwise limit of absolutely continuous solutions. 
  Let us point out that, in view of a result due to R.L. Baire (see \cite{Baire1899,Dun05} and also \cite[Remark 2.2]{AroRam14}), bounded measurable maps exist such that they are not the pointwise limit of continuous functions. An example of such a map is provided by the  {\em Dirichlet function ${\bf 1}_{\cQ\cap [0,1]}.$}
 }
\end{example}

\subsection{Consistency and uniqueness}
If $u$ is absolutely continuous, then for every control $v \in L^1([a,b];V),$ the corresponding Carath\'eodory solution is trivially a simple limit solution. A related question, namely whether the limit extension does or does not insert new absolutely continuous solutions,  is of obvious importance. By restricting attention to the first three variables in Example \ref{ex-suss}, one can notice that a simple limit solution may happen to be absolutely continuous, or even $\C^\infty$, without being a Carath\'eodory solution (a fact that turns out to be natural when translated in the language of {\it rough paths} \cite{Lyo94,Lyo02,Lej03}). 
However, {\it consistency} of limit solutions with Carath\'eodory solutions is achieved in two important cases: i) when the system is commutative and,  ii) when the limit solution is BV simple. Of course, this consistency implies uniqueness. More precisely:

\begin{theorem}
\label{cons} Let $\xb \in \cR^n,$ $(u,v)\in AC([a,b];U)\times L^1([a,b];V) $ and let us assume one of the following two conditions.
\begin{itemize}
\item[(A)] The system verifies Hypothesis (CC)  and $x$ is a limit solution corresponding to the initial condition $\xb$ and the control $(u,v).$
\item[(B)]
$x\in AC([a,b];\cR^n)$ is a BV simple limit  solution  corresponding to $\xb$ and $(u,v).$
\end{itemize}
Then $x$ coincides with $x[\xb,u,v],$ the (unique) Carath\'eodory solution corresponding to $\xb$ and $(u,v).$
\end{theorem}

Let us insist on the fact that assumptions in  (A) or (B) are crucial because of the noncommutativity of the vector fields $g_1,\dots,g_m.$
In fact, the statement corresponding to (A) is a straightforward consequence of a more general result on existence and uniqueness, valid for {\em all} limit solutions in the commutative case (see Theorem \ref{exunth}). While the statement  associated with condition (B) will be proved later, at the end of Section  \ref{noncommSec},
let us point out that the variations of the controls 
$u_k$ in Example \ref{ex-suss} tend to infinity as $k$ goes to infinity. So the fact that the first three components of $\tilde x$ are smooth does not contradicts statement (B) (when this is applied to the first three equations of system).



\subsection{Existence of  limit solutions}
Existence of limit solutions will be achieved as consequence of the representation results stated in the next sections. Let us begin with the commutative case.

 \begin{theorem}
 Let us assume that the commutativity hypothesis (CC) holds true. Then, for every $\xb\in \cR^n$ and  $(u,v)\in \L^1([a,b];U)\times L^1([a,b];V),$ there exists a limit solution of \eqref{E}-\eqref{IC}.
 \end{theorem}

 The result is contained in Theorem \ref{exunth} below.

 \vspace{5pt}

The noncommutative case is more involved. Yet, for a control with bounded variation, not only Theorem \ref{Existence2} below establishes the existence of a limit solution $x,$ but this $x$   turns out to be also a BV simple limit solution. To state this result we need the  definition of {\it Whitney property} for a subset of an Euclidian space, which relates the geodesic distance with the Euclidian distance.

\begin{definition} We say that a compact subset $U$ of $\cR^m$ {\em has the Whitney property} (see \cite{Whitney34})
 if there is $M\geq 1$ such that, for every pair  $(u_1,u_2)\in U\times U,$  there exists an absolutely continuous path  $ \gamma:[0,1]\to U$ verifying
$$
\gamma(0)=u_1,\quad\gamma(1)=u_2,\quad {\rm Var}_{[0,1]}(\gamma)\leq M |u_1-u_2|.
$$
\end{definition}

 A trivial example of a subset with the  Whitney property  is a compact star-shaped subset. Instead, the arcwise connected subset
$$
\{(x,y)\in\cR^2: x\in ]0,1], y=x\sin(1/x)\} \cup \{(0,0)\}
$$
does not have the  Whitney property.

\begin{theorem}
\label{Existence2}
Let us assume that $U$ has the Whitney property. Then, for every initial value $\xb \in \cR^n$ and control pair $(u,v) \in BV([a,b];U)\times L^1([a,b];V)$ there exists  an associated BV simple limit solution of \eqref{E}-\eqref{IC}.
 \end{theorem}

Theorem \ref{Existence2} will be proved in Section \ref{proofSec} as a consequence of  a representation property stated in Theorem \ref{gcpdgen1} below.

 Let us observe that, as a byproduct of Theorem \ref{Existence2}, we get the following density  result for BV functions in the Tychonoff topology of pointwise convergence.

 \begin{theorem}
 \label{whi}
 Let  us assume that $U$ has the Whitney property. Then, every map $u\in BV([a,b];U)$ is the pointwise limit of a sequence $(u_k)$ of absolutely continuous  maps $u_k:[a,b]\to U$  with equibounded variation.
 \end{theorem}

 Indeed, the latter result follows by Theorem \ref{Existence2} since the existence of a BV simple limit solution corresponding to $u$ implies the existence of a sequence with the mentioned properties.

\section {The commutative case } \label{commSec}

 In this section we mainly recall  some results from \cite{AroRam14}, where the commutativity   Hypothesis (CC) was assumed. In fact, as already stated in the Introduction, the notion of limit solution coincides with  a characterization of  the concept of  pointwise defined solution  formerly proved  for the particular case of  commutative systems \cite{BreRam91,AroRam14}. Let us
 also refer to  \cite{Sar91,Dyk94,Millerbook} for other references dealing with systems where the Lie algebra generated by $\{g_1,\dots,g_m\}$ is trivial.

\subsection{Existence, uniqueness, continuous dependence}

\begin{theorem}[Existence and Uniqueness]\label{exunth}
Assume that the Hypothesis (CC) holds.
Then,  for every initial value $\xb \in \cR^n$ and each control pair $(u,v)\in\L^1([a,b];U)\times L^1([a,b];V),$ there exists a unique limit  solution of the Cauchy problem \eqref{E}-\eqref{IC}.
\end{theorem}

Let us use $x[\xb,u,v]$ to denote the unique limit solution of \eqref{E}-\eqref{IC} corresponding to the initial value $\xb \in \cR^n,$ and a control pair $(u,v)\in \L^1([a,b];U)\times L^1([a,b];V).$

\begin{remark}
{\rm
In a more general situation when the state domain is a subset of $\cR^n,$ one can see that the only null bracket hypothesis $[ g_\alpha, g_\beta]\equiv 0$ does not guarantee the uniqueness of a limit solution (see an example in  \cite[Remark 1.1]{AroRam14}). 
}
\end{remark}

\begin{theorem}[Dependence on the data]
 \label{cont-depL1}
Let us assume Hypothesis (CC).  Then the following assertions hold true.
\begin{itemize}
\item[{(i)}] For each $\xb \in \cR^n$ and $u\in \L^1 ([a,b];U),$ the function $v\mapsto x[\xb,u,v]$ is continuous from $L^1([a,b];V)$ to $L^1([a,b];\cR^n).$
\item[(ii)] For any $r>0$ there exists a compact subset $K'\subset \cR^n,$ such that the trajectories $x[\bar x,u,v]$ have values in $K',$ whenever we consider $|\xb|\leq r,$  $u\in \L^1([a,b];U)$ and $v\in L^1([a,b];V).$
\item[(iii)] For each $r>0,$  there exists a constant $M>0$ such that, for every $t \in[a,b],$ for all $|\xb_1|, |\xb_2| \leq r,$ for all $u_1,u_2 \in \L^1([a,b];U)$  and for every $v \in L^1([a,b];V),$ one has
\benl
\begin{split}
&|x_1(t)-x_2(t)| +  \|x_1-x_2\|_1 \leq  \\
&\, M\Big[ |\xb_1-\xb_2|+ |u_1(a) - u_2(a)|  + |u_1(t)-u_2(t)|+\|u_1-u_2\|_1
\Big],
\end{split}
\eenl
where  $x_1 := x[\xb_1,u_1,v],$ $x_2 := x[\xb_2,u_2,v].$
\end{itemize}
\end{theorem}

\begin{remark}
Given that limit  solutions depend on the pointwise definition of $u,$ it is interesting to investigate the effects of a change of $u$ on a measure-zero subset of $[a,b].$ We refer to the reader to \cite[Theorem 2.4]{AroRam14} for this issue.
\end{remark}

 \section{The noncommutative case with BV controls}\label{noncommSec}

As soon as $[g_\alpha,g_\beta]\neq 0,$ the knowledge of a jump of the input $u$ is not enough to determine the jump of the  corresponding solution. 
 More precisely, when
$u$ has bounded variation:  at the discontinuity points  of $u,$  the non-drift dynamics and the  {\it $u$-path  during the jump} determine the corresponding discontinuity of the output. This fact is the common outcome of several investigations on the subject, which share a notion of solution  here referred  as  {\em graph completion solution}
 (see e.g. \cite{Ris65,BreRam88,DalRam91,MotRam95,SilVin96,CamFal99,PerSil00,Millerbook,Kar06,AruKarPer11}).  In general graph completion solutions are set-valued at a countable subset of instants and here are referred as  {\em set-valued graph completion solutions}, while
 their selection will be called simply {\em graph completion solution.}  The main result of this section is Theorem \ref{gcpdgen1} which establishes that the concepts of graph completion solution and that of BV simple limit solution are in fact equivalent.

 \subsection{Graph completion solutions}
 For any $L>0,$ we  use $\U_L$  to denote the subset of $L$-Lipschitz  maps $(\ubf_0,\ubf)$ belonging to ${\rm Lip}( [0,1];[a,b]\times U)$ such that $\ubf_0:[0,1]\to [a,b]$  is increasing and surjective, and we set $\U:= \bigcup_{L\geq 0}\U_L.$ 

\begin{definition}[Space-time control]
The elements $(\varphi_0,\varphi,\psi)$ of the set $\U\times L^1([0,1];V)$ will be called  {\em space-time controls}.
\end{definition}

Let us consider the {\it space-time control system} in the interval $[0,1]$ defined as
\bel{spacetime}
\left\{
\begin{split}
\ybf_0'(s) &= \ubf_0'(s),\\
\ybf'(s) &= f(y_0(s),\ybf(s),\ubf(s),\vbf(s)) \ubf_0'(s)+ \sum_{\alpha=1}^m {g}_\alpha(\ybf(s)) {\ubf_\alpha}'(s) , \\
(\ybf_0,\ybf&)(0) = (a,\bar x)\,,
\end{split}
\right.
\ee
where the apex $'$ denotes differentiation with respect to the  {\it pseudo-time} $s\in [0,1],$  and the (space-time) controls $(\ubf_0,\ubf,\vbf )$ belong to $\U\times L^1([0,1];V).$ For the sake of simplicity, we shall fix the initial value $\bar x \in \cR^{n}.$
For a given initial value $\xb\in\cR^n,$ let $\I$ denote the input-output map associated with the system \eqref{spacetime}, which assigns to each $(\varphi_0,\varphi,\psi)\in \U \times L^1([0,1];V)$ verifying $(\varphi_0,\varphi)(0)=(a,\ub),$
 the corresponding solution $(y_0,y)[\varphi_0,\varphi,\psi]$ of the Cauchy problem \eqref{spacetime}. The following property holds:

\begin{proposition}
\label{Icont}
For every $L>0,$ the map $\I$ is continuous from $\U_L\times L^1([0,1];V)$ to ${\rm Lip}([0,1];[a,b]\times \cR^n),$ where the domain and the range are endowed with the product topology of $\C^0([0,1];[a,b]\times U) \times L^1([0,1];V)$ and the $\C^0([0,1];[a,b]\times \cR^n)-$topology, respectively.
\end{proposition}

The previous result is an extension of \cite[Theorem 1]{BreRam88} (where $f$ is independent of $v$). For a proof we refer to  the first part of the proof of  \cite[Theorem 4.1]{MotRam95}.

\if{ 
 \begin{remark}\label {paramfree}
 {\rm
 Let us point out that the space-time system \eqref{spacetime}
has {\it parameter free} character. We mean that  if $\hat s:[0,1]\to [0,1]$ is a bijective, increasing, bi-Lipschitz map then one has
$$
(y_0,y)[\ubf_0,\ubf,\psi]\circ \hat s = (y_0,y)[\ubf_0\circ \hat s,\ubf\circ \hat s,\psi\circ \hat s].
$$
}\end{remark}
}\fi

Consider now $(u,v) \in AC([a,b];U)\times L^1 ([a,b];V),$ and $\ubf_0:[0,1]\to[a,b]$  a Lipschitz, bijective map such that $\varphi_0(0)=a$ and $\ubf:= u\circ \varphi_0$ is Lipschitz continuous.
It is not restrictive to assume that $v$ is Borel measurable, since in every $L^1$ equivalence class there are always Borel measurable representatives.
Then, with $v$ being $L^1$ and Borel measurable, the composition $\psi:= v \circ \varphi_0$ belongs to $L^1([0,1];V);$ and the unique corresponding solution  $(y_0,y)[\varphi_0,\varphi,\psi]$ of \eqref{spacetime} verifies
$$
y(s) = x\circ \varphi_0(s), \quad \forall s\in [0,1],
$$
where $x$ is the Carath\'eodory solution of \eqref{E}-\eqref{IC} corresponding to $(u,v).$
By means of the notion of {\it graph completion}, this argument can be  used to extend the notion of  output of the system \eqref{E}  to controls $u$ with bounded variation (and possibly discontinuous).

\begin{definition}[Graph completion \cite{BreRam88}]
Consider a control $u\in BV([a,b];U).$ A {\em graph completion} (shortly {\em g.c.}) of $u$ is  a space-time control $(\ubf_0,\ubf)\in \U$ such that,
$\forall t\in [a,b],$ there exists $s\in [0,1]$ verifying $(t,u(t)) = (\ubf_0,\ubf)(s).$
\end{definition}

Consider $(u,v) \in BV([a,b];U) \times L^1 ([a,b];V),$  let $(\varphi_0,\varphi)\in \U$ be a g.c. of $u$ and set $\psi:=v \circ \varphi_0.$ Then either $(y_0,y)[\varphi_0,\varphi,\psi]$ or $(y_0,y)[\varphi_0,\varphi;v]$ shall denote the solution of \eqref{spacetime} corresponding to $(\varphi_0,\varphi,\psi).$

\begin{definition}[Set-valued graph completion solution \cite{MotRam95}]
\label{gcs}
Consider $(\ubf_0,\ubf)$ a g.c. of $u\in BV([a,b];U),$ $v$ a control in $L^1([a,b];V),$ and $(y_0,y)[\varphi_0,\varphi;v]$  the corresponding solution of \eqref{spacetime}. We say that the (possibly set-valued) map  defined by
\benl
x(t):= \ybf\circ\ubf_0^{\leftarrow}(t),\qquad \forall t\in [a,b],
\eenl
is the\footnote{Uniqueness is a consequence of the uniqueness of solutions to Cauchy problem \eqref{spacetime}. } {\em set-valued graph completion  solution} to \eqref{E}-\eqref{IC}, shortly {\em set-valued g.c. solution,} corresponding to the control triple $(\ubf_0,\ubf;v).$ Here $\ubf_0^{\leftarrow}(t)$ denotes the pre-image of the singleton $\{t\}.$
When the initial condition is meant by the context this solution will be denoted by $x[\ubf_0,\ubf;v].$
\end{definition}

\begin{remark} {\rm A set-valued g.c. solution $x=x[\ubf_0,\ubf;v]$ is single-valued at almost every $t,$ for there are at most countably many intervals $I\subset [0,1]$ on which $\ubf_0 $ is constant and $\ubf$ is not constant (this being a consequence of $(\ubf_0,\ubf)\in BV$)\footnote{Let us point out that $x$ may well be set-valued even at times $t$ where $u$ is continuous: unless the system is globally integrable  for this to happen, it is necessary that $\ubf_0^\leftarrow(t) = [s_1,s_2], s_1<s_2,$ and $\ubf_{\big|_{[s_1,s_2]}}$ is a non trivial closed curve.}. 
}
\end{remark}


\subsection{Equivalence between g.c. solutions and BV simple limit solutions}
Let us   give the definition of {g.c. solution} by means of a selection of the inverse (set-valued) map $t\mapsto (\ubf_0,\ubf)^{\leftarrow}(t).$   We call such a selection a {\it $(\ubf_0,\ubf)$-clock.}

\begin{definition} 
\label{DefClock}
Let us consider  $u\in BV([a,b];U)$ and let $(\ubf_0,\ubf)$ be a graph completion of $u.$ Any map   $\sigma: [a,b]\to [0,1]$ such that, $\forall t\in[a,b],$
$$
(\ubf_0,\ubf)\circ\sigma(t) = (t,u(t))
$$
will be called a  {\em clock}  corresponding to  $(\ubf_0,\ubf).$
\end{definition}
\begin{remark}
\label{clockremark}
{\rm
Notice that if $L>0$ is a Lipschitz  constant for $\ubf_0,$ then any clock $\sigma$ verifies
$$
\sigma(t_2)-\sigma(t_1)\geq \frac{1}{L}(t_2-t_1),
$$
for all $t_1,t_2\in [a,b]$ such that $t_1 < t_2.$
}
\end{remark}

\begin{definition}
\label{ugcs}
Let $(\ubf_0,\ubf)$ be a g.c. of a control $u$ in $BV ([a,b];U),$ and let  ${v} \in L^1([a,b];V).$ Let $\sigma$ be a clock corresponding to $(\ubf_0,\ubf).$
The map
\benl
x(t):= \ybf[\ubf_0,\varphi;v]\circ \sigma(t),\qquad \forall t\in [a,b],
\eenl
will be called  the   {\em g.c. solution} corresponding to $(\ubf_0,\ubf;v)$ and to the clock $\sigma.$
We shall use $x[\ubf_0,\ubf;v]_{ \sigma}$ to denote this solution.
\end{definition}

The notion of g.c. solution  in fact coincides with that  of   limit solution, as stated in the following theorem.

\begin{theorem}
\label{gcpdgen1}
Let $\xb \in \cR^n$ and $(u,v)\in BV([a,b];U) \times L^1([a,b];V).$
A map $x:[a,b]\to\cR^n$ is a g.c. solution of \eqref{E}-\eqref{IC} corresponding to $(\xb,u,v)$
if and only if $x$ is a BV simple limit solution of \eqref{E}-\eqref{IC}.
\end{theorem}

We postpone the proof of this theorem to  the next section.

\begin{remark}
{\rm
A characterization of graph completion solutions was also provided in \cite{SilVin96}, where it is shown  equivalence of the latter  with {\it robust solutions} {for a class of differential inclusions with scalar impulsive control}.
}
\end{remark}

\begin{remark} {\rm
 The ``only if" part of Theorem \ref{gcpdgen1} differs from former  approximation results of graph completion solutions (see e.g. \cite{WolZab07} for a careful discussion on the subject, in particular, see Theorem 4.1 therein) in that here we consider  a pointwise (hence non-metric) convergence. This is obviously due to the fact that our notion of graph completion solution is  defined at {\it all } points of the interval $[a,b]$ (customarily only left and right limits are considered at jump points). In turn,  on one hand this comes from the aim of having a unified notion of solution, matching with the (everywhere defined) concept of solution of the commutative case.  On the other hand graph completions allow for jumps of the trajectory even at the times $t$  where $u$ is continuous (a {\it loop} of $u$ could be considered at these instants, which, thanks to the non-triviality of the Lie algebra generated by $\{g_1,\dots,g_m\}$,  might  generate a discontinuity in $x$). 

 Finally notice that the  ``if" part is not an approximation result but rather an existence result, still rendered not obvious by the fact that we are considering pointwise convergence at any $t\in [a,b]$.}
\end{remark}

\subsection{ Proof of Theorem \ref{cons}. }

Let us conclude this section with the proof of the consistency result in the noncommutative case, namely the case (B) of  Theorem \ref{cons}.
Let us recall its statement:

\noindent  {\em Let $\xb \in \cR^n,$ $(u,v)\in AC([a,b];U)\times L^1([a,b];V) $ and let us assume that $x\in AC([a,b];\cR^n)$ is a BV simple limit  solution  corresponding to $\xb$ and $(u,v).$
Then $x$ coincides with $x[\xb,u,v],$ the (unique) Carath\'eodory solution corresponding to $\xb$ and $(u,v).$}

\begin{proof}
In view of the ``if'' part of Theorem \ref{gcpdgen1}, we know that $x$ is a g.c. solution. Let $(\varphi_0,\hat\varphi)$ be a g.c. of $u$ and $\sigma$ be a clock such that
\be
\label{gcformula}
(t,u(t),x(t)) = (\varphi_0,\hat\varphi,\hat y)\circ \sigma (t),\quad \text{for all } t\in[a,b],
\ee
where $\hat y:= y[\varphi_0,\hat\varphi;v].$

Set now $\varphi:=u\circ  \varphi_0,$ and consider the g.c. $(\varphi_0,\varphi).$ We shall prove next that the g.c. solution associated to $(\varphi_0,\varphi)$ and the clock $\sigma$ coincides with the Carath\'eodory solution $\tilde x:= x[\xb,u,v]$ of \eqref{E}-\eqref{IC}. Finally, we will see that $x=\tilde x,$ from where the desired result will follow.

Let $\tilde\varphi_{0,k}:[0,1]\to [a,b]$ be a sequence of Lipschitz continuous strictly increasing functions such that
$$
\tilde\varphi_{0,k} \to \varphi_0\quad \text{uniformly, } \,\, \tilde\sigma_k:=\tilde\varphi_{0,k}^{-1} \to \sigma\quad \text{pointwise.}
$$
The existence of such a sequence $(\tilde\varphi_{0,k})$ is guaranteed by Proposition \ref{Propsigma}.
Set $\tilde\varphi_k:=u\circ \tilde\varphi_{0,k}.$ Since $u$ is uniformly continuous on $[a,b],$ one has
$$
\tilde\varphi_k = u\circ \tilde\varphi_{0,k} \to \varphi,\quad \text{uniformly.}
$$
Hence, in view of Proposition \ref{Icont}, the corresponding solutions $\tilde y_k:= y[\tilde\varphi_{0,k},\tilde\varphi_k;v] $ converge uniformly to $y[\varphi_0,\varphi;v].$ But, for all $k,$ $\tilde y_k \circ \tilde\sigma_k=\tilde x$ and, therefore, one gets
$$
\tilde x(t) = \tilde y_k \circ \tilde \sigma_k(t) \to y \circ \sigma (t), \quad \text{for all  } t\in [a,b],
$$
due to the uniform convergence of $\tilde y_k$ and since $\tilde \sigma_k$ tends to $\sigma$ pointwise.
Therefore, the g.c. solution $y \circ \sigma$ is equal to $\tilde x.$

Let us go back to the g.c. $(\varphi_0,\hat\varphi).$ We know that $\varphi_0$ is increasing on $[0,1]$ and constant on the intervals $[s_{1,i}:=\sigma(t_i-),s_{2,i}:=\sigma(t_i+)],$ where $t_i$ are the (at most countable) discontinuity points of $\sigma.$
Notice that
\be
\label{varphisi}
(\hat\varphi,\hat y)(s_{1,i})= \lim_{t \nearrow t_i} (\hat\varphi,\hat y)\circ \sigma(t) = \lim_{t \nearrow t_i}  (u,x)(t) = (u,x)(t_i),
\ee
where the second equality follows from \eqref{gcformula}, and the third one holds due to the continuity of $(u,x).$ Analogously,
\be
\label{varphisi2}
(\hat\varphi,\hat y)(s_{2,i})=(u,x)(t_i).
\ee
Thus, for all $i\in \cN,$
$$
(t_i,u(t_i),x(t_i)) = (\varphi_0,\hat\varphi,\hat y)( s_{1,i}) = (\varphi_0,\hat\varphi,\hat y)\circ \sigma (t_i-),
$$
and, therefore, one can redefine $\sigma$ to be left-continuous and still satisfy \eqref{gcformula} (i.e. the resulting $\sigma$ is a $(\varphi_0,\hat\varphi)-$clock).

In the sequel, we prove that the g.c. solution $x$ coincides with the Carath\'eodory solution $\tilde x.$
Set $\varphi_1=\hat\varphi$ and, for $k\in \cN,$ define (recursively)
\be
\varphi_{k+1}(s):=
\left\{
\ba{cl}
\varphi(s), &\text{ if } s\in [s_{1,k},s_{2,k}],\\
\varphi_{k}(s), &\text{ otherwise}.
\ea
\right.
\ee
Note that each $\varphi_k$ is continuous in view of \eqref{varphisi}-\eqref{varphisi2}. Actually, each $(\varphi_0,\varphi_k)$ is Lipschitz continuous (with constant bounded by the Lipschitz constant of $(\varphi_0,\hat\varphi)$) and is a g.c. of $u.$
Moreover, $\varphi_k $ converges pointwise to $\varphi$ on $[0,1]$ and, for $j>k,$ one has
$$
\|\varphi_{j} - \varphi_k\|_{1,1} = \sum_{i=k}^{j-1} \int_{s_{1,i}}^{s_{2,i}} \big( |\varphi_{i+1} - \varphi_i| + |\varphi_{i+1}' - \varphi_i'| \big) {\rm d}s \leq  \sum_{i=k+1}^j (s_{2,i}-s_{1,i})M,
$$
where $M>0$ is such that $\|\hat\varphi\|_\infty+\|\hat\varphi'\|_\infty \leq M,$ and $\|\cdot\|_{1,1}$ denotes the norm in the Sobolev space $  W^{1,1}([0,1];\cR^m).$
Then one gets
$$
\varphi_k \to \varphi \ \text{in } W^{1,1}([0,1];\cR^m),
$$
since
$$
\sum_{i=k}^\infty (s_{2,k}-s_{1,k}) \to 0, \quad \text{when $k\to\infty$}.
$$
Hence, the corresponding solutions
$
y_k:=y[\varphi_0,\varphi_k;v]
$
converge to $y$ uniformly. Let us show that all the $y_k \circ \sigma$ are equal to $x.$ By definition $y_1 \circ \sigma = \hat y \circ \sigma = x.$ Now, to pass to $y_2,$ one modifies $\varphi_1$ on the interval $[s_{1,1},s_{2,1}].$ Both $\varphi_0$ and the control $\varphi_2$ are constant on $[s_{1,1},s_{2,1}].$ Hence $y_2(s_{2,1})=y_2(s_{1,1})=x(t_1).$ Thus, $y_1$ and $y_2$ may differ only on $]s_{1,1},s_{2,1}[.$  Hence, $ y_2 \circ \sigma = y_1 \circ \sigma = \hat y \circ \sigma,$ i.e. the g.c. solution associated to $\varphi_1=\hat\varphi$ (and the clock $\sigma$) coincides with the g.c. solution associated to $\varphi_2.$  Analogously, one can prove that $y_k \circ \sigma = \hat y \circ \sigma$ for all $k,$ hence all the g.c. solutions coincide. Thus,
$$
x = \hat y \circ \sigma = y_k \circ \sigma \to y \circ \sigma = \tilde x,
$$
and the result follows.
 \end{proof}

\section{Proofs of Theorems \ref{gcpdgen1} and \ref{Existence2}}\label{proofSec}

\subsection{Proof of the ``only if'' part of Theorem \ref{gcpdgen1}}

Let us begin by proving the following result on pointwise approximation of increasing maps:

\begin{theorem}
\label{Propsigma}
Let $\sigma:[a,b]\to [0,1]$ be a strictly increasing function such that
\be
\label{sigmaL}
\sigma(t_2)-\sigma(t_1) \geq \frac{1}{L}(t_2-t_1), \text{ for every pair $t_1,t_2 \in [a,b],$ with $t_1 < t_2,$}
\ee
 for some $L>0.$ Then, there exists a sequence of  functions $(\sigma_k)$ with the following properties:
\begin{itemize}
\item[(i)] for each $k\in \cN,$  $\sigma_k$ is a strictly increasing, absolutely continuous function from $[a,b]$ onto $[0,1];$
\item[(ii)]
\be
\label{limitsigmak}
|\sigma_k(t) - \sigma(t)| + \|\sigma_k-\sigma\|_1 \to 0,\quad
\text{for all } t\in [a,b];
\ee
\item[(iii)] the inverses $\varphi_{0,k}:= \sigma_k^{-1}$ (are strictly increasing, $L-$Lipschitz continuous and) converge uniformly to the unique ($L-$Lipschitz continuous) increasing surjective map $\varphi_0:[0,1]\to [a,b],$ verifying
$$
\varphi_0 \circ \sigma(t) =t,\quad \forall t\in [a,b].
$$
\end{itemize}
\end{theorem}

{\em Proof.}
Let us start by proving the following claim:
\begin{claim}
\label{ClaimMollifier}
There exists a sequence $(\hat\sigma_k)$ of strictly increasing functions from $[a,b]$ onto $[0,1]$ with the following properties:
\begin{itemize}
\item[(a)] each $\hat\sigma_k$ is $\C^\infty,$  and verifies
\be
\label{derhatsigma}
\frac{{\rm d}\hat\sigma_k(t)}{{\rm d}t} \geq \frac{1}{L},\quad \forall t\in [a,b];
\ee
\item[(b)] for all $t\in [a,b],$ one has
\be
\label{limhat}
\left|\hat\sigma_k(t) - \frac{\sigma(t-)+\sigma(t+)}{2}\right|
+ \|\hat\sigma_k - \sigma\|_1 \to 0;
\ee
\item[(c)] the inverses $\hat\varphi_{0,k}:= \hat\sigma_k^{-1}$ (are strictly increasing, $L-$Lipschitz continuous and) converge uniformly to $\varphi_0.$
\end{itemize}
\end{claim}

{\em Proof of Claim \ref{ClaimMollifier}.}
Let us choose $M\in ]0,b-a]$ and let us first extend $\sigma$
to $[a-M,b+M]$ by setting  $\sigma(a-t):=-\sigma(a+t),$ and $\sigma(b+t) : =- \sigma(b-t),$ for $t\in ]0,M].$
Define $\sigma$ on the whole real line $\cR$ by setting $\sigma(t)=0,$ for $t\in
\cR \backslash [a-M,b+M].$

Let $\rho:\cR\to [0,+\infty[$ be a $\C^{\infty},$ even mapping such that ${\rm supp} (\rho) \subseteq [-M,M]$ and $\|\rho \|_1=1.$   For every $k\in \cN,$  set $\rho_k(t):= 2k\,\rho (2kt),$ for all $t \in \cR.$ Define  the approximating map  $\hat\sigma_k:[a,b]\to\cR$
by setting
\be
\label{sigmadef}
\hat\sigma_k(t):= (\sigma \ast\rho_k)(t):=\int_{-\infty}^\infty \sigma(t-\tau)\rho_k(\tau){\rm d}\tau,
\quad \forall t\in[a,b].
\ee
The $\C^\infty$ regularity of $\hat\sigma_k$  claimed in (i) is a well-known result and its strict monotonicity is just a consequence of the strict monotonicity of ${\sigma}$  (for $\rho_k$ is nonnegative).
Since the mapping $t\mapsto \sigma(a-t)$ is odd and $t\mapsto \rho_k(t)$ is even, one gets
$$
\hat\sigma_k(a) = \int_{-\infty}^{\infty}  \sigma(a-t)\rho_k(t) dt = 0,
$$
Similarly,
$$
\hat\sigma_k(b) = 1.
$$

Let us show that \eqref{derhatsigma} holds true. Indeed, if $t_1<t_2$ are points in $[a,b],$ one has
\benl
\begin{split}
\hat\sigma_k(t_2) - \hat\sigma_k(t_1) &= \int_{-\infty}^\infty \big(\sigma(t_2-t) - \sigma(t_1-t) \big) \rho_k(t) {\rm d} t \\
& \geq \int_{-\infty}^\infty \frac{1}{L} (t_2-t_1)\rho_k(t) {\rm d} t = \frac{1}{L} (t_2-t_1),
\end{split}
\eenl
where  the first inequality follows from \eqref{sigmaL}.

Let us prove (b). For any $t \in [a,b]$ and $k\in \cN,$  one has
\be
\label{limhatsigma}
\begin{split}
&\hat\sigma_k(t) - \frac{\sigma(t-)+\sigma(t+)}{2}\\
&= \int_{-M/2k}^{0} \sigma(t -\tau)\rho_k(\tau) {\rm d}\tau
+\int_0^{M/2k} \sigma(t -\tau)\rho_k(\tau) {\rm d}\tau - \frac{\sigma(t-)+\sigma(t+)}{2}\\
&= \int_{-M}^{0} [\sigma(t -\theta/2k)-\sigma(t+)]\rho(\theta) {\rm d}\theta
+\int_0^{M} [\sigma(t -\theta/2k)-\sigma(t-)]\rho(\theta) {\rm d} \theta.
\end{split}
\ee
The functions
$$
\theta\mapsto \sigma(t -\theta/2k)-\sigma(t+)
$$
are equibounded by $2,$ and $\lim_{k\to \infty} \sigma(t -\theta/2k)-\sigma(t+)=0,$ for all $\theta \in [-M,0[.$
Analogous facts hold if we replace $t+$ with $t-.$
From \eqref{limhatsigma} and the Dominated Convergence Theorem, we get that
$$
\left|\hat\sigma_k(t) - \frac{\sigma(t-)+\sigma(t+)}{2}\right| \to 0,\quad \text{for all } t\in [a,b].
$$
The latter limit implies that $\hat\sigma_k(t) \to \sigma(t)$ at each $t$ where $\sigma$ is continuous, hence, almost everywhere on $[a,b].$
Since $|\hat\sigma_k-\sigma|$ is bounded by 2 on $[a,b],$ once again by the Dominated Convergence Theorem, we get that $\|\hat\sigma_k-\sigma\|_1 \to 0,$ so (b) is proved.

As for (c), note that the inverses $\hat\varphi_{0,k}:=\hat\sigma_k^{-1}$ are  strictly increasing, $L-$Lipschitz continuous and that, taking if necessary a subsequence, they converge to an $L-$Lipschitz continuous function $\hat\varphi_0.$ We shall prove that $\hat\varphi_0 \circ \sigma (t)=t,$ on $[a,b],$ which implies that $\hat\varphi_0 = \varphi_0.$
Indeed, observe that, due to (b),   for any continuity point $t\in [a,b]$ of $\sigma$ one has
\be
\label{thetasigma}
 t= \hat\varphi_{0,k} \circ \hat\sigma_{k} (t)  \to \hat\varphi_0(\sigma(t)).
\ee
If $t$ is not  a continuity point of $\sigma,$ let $(t^1_k)$ and $(t^2_k)$ be two sequences in $[a,b]$ such that $t^i_k$ is a continuity point of $\sigma$ for $i=1,2$ and $k\in \cN,$ and
\benl
t^1_k < t < t^2_k,\quad t_k^1 \to t,\quad t_k^2\to t .
\eenl
Given that  both $\sigma$  and $\hat\varphi_0$ are increasing, one has, in particular,
$$
\hat\varphi_0(\sigma(t_k^1)) \leq \hat\varphi_0 ( \sigma(t)) \leq \hat\varphi_0 ( \sigma(t^2_k)).
$$
Thus, in view of \eqref{thetasigma}, one obtains
$$
t_k^1 \leq \hat\varphi_0(\sigma(t)) \leq t_k^2, \quad \text{for all } k\in \cN.
$$
Taking the limit as $k$ goes to infinity we conclude that $\hat\varphi_0 (\sigma(t) ) = t,$  as desired.
This completes the proof of the Claim \ref{ClaimMollifier}.

\vspace{5pt}

Let $\T \subset [a,b]$ be the set of discontinuity points of $\sigma. $ Since $\T$ is at most countable, we can write $\T = \{ t_i: i\in \cN\}.$ Let $t_i \in \T,$ and set $s_{1,i} :=\sigma(t_i-),$ and $s_{2,i} := \sigma(t_i+).$ Since $\sigma$ is  increasing, one gets
$$
s_{1,i} \leq  \sigma(t_i)  \leq s_{2,i}.
$$
According to Claim \ref{ClaimMollifier}, one has
$
\hat\sigma_k(t_i) \to \ds\frac{s_{1,i}+s_{2,i}}{2},
$
as $k\in \infty.$
Hence, there exists $ k_i \in \cN,$ such that,
for $k \geq k_i,$ the inequality $s_{1,i} < \hat\sigma_k(t_i) < s_{2,i}$ holds and, consequently,
$$
\hat\varphi_{0,k}(s_{1,i}) < t_i < \hat\varphi_{0,k}(s_{2,i}).
$$
If one has
\be
\label{sigmainterior}
s_{1,i} <  \sigma(t_i)  < s_{2,i},
\ee
define, for  $k \geq k_i$ the function $\tilde{\varphi}_{0,k}$ on $ [s_{1,i},s_{2,i}]$ by letting
\benl
\tilde{\varphi}_{0,k}(s):=
\left\{
\ba{cl}
\hat\varphi_{0,k}(s_{1,i}) + \ds\frac{t_i-\hat\varphi_{0,k}(s_{1,i})}{\sigma(t_i)-s_{1,i}} (s-s_{1,i}),\quad & \text{for }s\in [s_{1,i},\sigma(t_i)],\\
t_i +\ds \frac{\hat\varphi_{0,k}(s_{2,i})-t_i}{s_{2,i}-\sigma(t_i)}(s-\sigma(t_i)),\quad
& \text{for } s\in ]\sigma(t_i),s_{2,i}].
\ea
\right.
\eenl
In particular, one gets
\be
\label{sigmati1}
\tilde \varphi_{0,k}^{-1}(t_i)=\sigma(t_i).
\ee
Note that, on $[s_{1,i},s_{2,i}],$ the function $\varphi_0$ is constantly equal to $t_i.$
Since $\tilde\varphi_{0,k}(s_{1,i}) \to \varphi_0(s_{1,i})= t_i,$ $\tilde\varphi_{0,k}(s_2,i)\to \varphi_0(s_{2,i})= t_i,$ and $\tilde\varphi_{0,k}$ is increasing and piecewise affine,  there exists $\tilde k_i \geq k_i$ such that
\be\label{dtilde}
\left\| \frac{d}{ds}{\tilde{\varphi}_{0,k}}(s) \right\|_\infty \leq L,\quad \text{on } [s_{1,i},s_{2,i}],\ \forall k \geq \tilde k_i.
\ee
If, instead of \eqref{sigmainterior}, one has that $\sigma(t_i) = s_{1,i},$ then set, for $k \geq k_i,$
\benl
\tilde{\varphi}_{0,k}(s):=
\left\{
\ba{cl}
\hat\varphi_{0,k}(s_{1,i}) + L (s-s_{1,i}),\quad & \text{for } s\in [s_{1,i},\hat s_{i,k}],\\
t_i +\ds \frac{\hat\varphi_{0,k}(s_{2,i})-t_i}{s_{2,i}-\hat s_{i,k}}(s-\hat s_{i,k}),\quad
& \text{for } s\in ]\hat s_{i,k},s_{2,i}],
\ea
\right.
\eenl
where $\hat s_{i,k}\in ]s_{1,i},s_{2,i}]$ is determined by $ \hat\varphi_{0,k}(s_{1,i}) + L (\hat s_{i,k}-s_{1,i})=t_i.$  This $\hat s_{i,k}$ exists, for $\hat\ubf_{0,k}$ is $L$-Lipschitz continuous.
Notice that
\be
\label{sigmati2}
\tilde\varphi_{0,k}^{-1}(t_i) \to s_{1,i}=\sigma(t_i), \text{ when } k\to \infty.
\ee
Furthermore, for some $\tilde k_i \geq k_i,$ one gets  \eqref{dtilde} as well . An analogous argument can be applied when $\sigma(t_i) = s_{2,i}.$

Let us construct a sequence $(\varphi_{0,k})$ of maps from $[0,1]$ onto $[a,b]$ by setting,
if $s\in [s_{1,i},s_{2,i}]$ for some $i,$
\be
\varphi_{0,k}(s):=
\left\{
\ba{cl}
\tilde{\varphi}_{0,k}(s),\quad &\text{if } k \geq \tilde k_i,\\
\hat{\varphi}_{0,k}(s),\quad &\text{otherwise},
\ea
\right.
\ee
and letting  $\varphi_{0,k}(s):=\hat\varphi_{0,k}(s) $ for all $s \in [0,1]\backslash \bigcup_{i=1}^\infty [s_{1,i},s_{2,i}].$

It is easy to verify that the functions $\varphi_{0,k}$ are equi-Lipschitz continuous (with Lipschitz constant bounded by $L$), strictly increasing and surjective. Moreover,  they converge uniformly to $\varphi_0.$ Then (i) and (iii) hold true if one sets  $\sigma_k:=\varphi_{0,k}^{-1}.$

To prove the convergence stated in  (ii), begin by observing that, in view of \eqref{sigmati1} and \eqref{sigmati2}, at each  $t_i\in \T$ one has
\be
\label{convdisc}
\sigma_k(t_i) \to \sigma(t_i).
\ee
Suppose now that $t \in ]a,b[ \backslash \T.$ Notice that if $t$ is such that
\be
\label{tnotlimit}
t \in ]a,b[ \backslash \bigcup_{i=1}^\infty \varphi_{0,k}( [s_{1,i},s_{2,i}]), \text{ for all $k$ sufficiently large},
\ee
then
$$
\sigma_k(t) = \varphi_{0,k}^{-1}(t) = \hat\varphi_{0,k}^{-1}(t) = \hat\sigma_k(t).
$$
Thus, from the pointwise convergence $\hat\sigma_k(t) \to \sigma(t)$ in Claim \ref{ClaimMollifier}, one deduces that $\sigma_k(t) \to \sigma(t).$
If, on the contrary, \eqref{tnotlimit} does not hold, then for every $\tilde k,$ there exists $k \geq \tilde k$ and $i_k \in \cN,$ such that
\be
\label{tlimit}
\varphi_{0,k}( s_{1,i_k}) \leq t\leq \varphi_{0,k}( s_{2,i_k}).
\ee
We can assume, without loss of generality, that \eqref{tlimit} holds true for all $k\in \cN,$ and for a sequence $(i_k).$ Notice that, in view of \eqref{tlimit}, one has
\benl
\begin{split}
|t_{i_k}-t| \leq &|t_{i_k}-\varphi_{0,k}(s_{1,i_k})|+|t_{i_k}-\varphi_{0,k}(s_{2,i_k})|
=  \\
&|\varphi_{0}(s_{1,i_k})-\varphi_{0,k}(s_{1,i_k})|+|\varphi_{0}(s_{2,i_k})-\varphi_{0,k}(s_{2,i_k})|.
\end{split}
\eenl
Thus,
\be
\label{limth}
t_{i_k} \to t,
\ee
since $\varphi_{0,k}\to \varphi_0$ uniformly.

By definition of $s_{1,i_k},$ there exists, for each $k,$ $\eps_{i_k}>0$ such that for all $\theta \in ] t_{i_k}-\eps_{i_k}, t_{i_k}[,$ one has
$|\sigma(\theta)-s_{1,i_k}| < 1/2k.$ Thus, for each $k,$ there exists $\theta_{i_k} \in ]a, t_{i_k}[$ verifying
\be
\label{thetaik}
|\theta_{i_k} - t_{i_k}|+|\sigma(\theta_{i_k})-s_{1,i_k}| < 1/k.
\ee
Consequently, from latter equation and \eqref{limth}, one obtains
$\theta_{i_k} \to t.$ Moreover, the continuity of $\sigma$ at $t,$ implies that $\sigma(\theta_{i_k}) \to \sigma(t),$ and, in view of \eqref{thetaik}, one has $s_{1,i_k} \to \sigma(t).$ Analogously, one can prove that $s_{2,i_k} \to \sigma(t).$ Finally, observe that \eqref{tlimit} holds if and only if
$$
s_{1,i_k} \leq \sigma_k(t) \leq s_{2,i_k},
$$
and, since both the right and left hand-sides of previous inequality converge to $\sigma(t),$ one concludes that $\sigma_k(t) \to \sigma(t)$ as well.
\if{
then for some $i\in \cN$ and $k$ sufficiently large, one has
$$
\varphi_{0,k}(s_{2,i}) < t < \varphi_{0,k}(s_{1,i+1}).
$$
Since the functions $\hat\varphi_{0,k}$ and $\varphi_{0,k}$ coincide on $[s_{2,i},s_{1,i+1}],$  one obtains
$$
\sigma_k(t) = \varphi_{0,k}^{-1}(t) = \hat\varphi_{0,k}^{-1}(t) = \hat\sigma_k(t),\quad \text{for $k$ large enough.}
$$
Thus, from the pointwise convergence $\hat\sigma_k(t) \to \sigma(t)$ in Claim \ref{ClaimMollifier}, one deduces that $\sigma_k(t) \to \sigma(t),$ for all $t\in \T.$
}\fi
The latter assertion together with equation \eqref{convdisc} and the Dominated Convergence Theorem yield (ii).
This completes the proof of Theorem \ref{Propsigma}.
\epf


Let us now prove the ``only if'' part of Theorem \ref{gcpdgen1}. Namely we shall show the following fact: 

{\em If $x:[a,b]\to\cR^n$ is a g.c. solution of \eqref{E}-\eqref{IC} associated to $\xb \in \cR^n$ and $(u,v)\in BV([a,b];U) \times L^1([a,b];V),$ then $x$ is a BV simple limit solution.}

\vspace{5pt}

Let $x:[a,b] \to \Omega$ be a g.c. solution corresponding to $(u,v).$ This means that there exists a graph completion $(\varphi_0,\varphi)$ of $u,$ and a clock $\sigma$ such that
$$
(t,u(t),x(t))= (\varphi_0,\varphi,y) \circ \sigma(t),\quad \forall t\in [a,b],
$$
where  $y:=y[\varphi_0,\varphi;v].$
For this clock $\sigma,$ consider the associated sequences $(\sigma_k)$ and $(\varphi_{0,k})$ provided by Theorem \ref{Propsigma}. Then the $(\varphi_{0,k},\varphi)$ are equi-Lipschitz continuous and
converge uniformly to $(\varphi_0,\varphi),$ so the corresponding solutions $y_k:=y[\varphi_{0,k},\varphi;v]  $ converge uniformly to $y$ (see \cite[Theorem 1]{BreRam88}).
Set, for any $t\in [a,b],$
$$
u_k(t) : = \varphi\circ \sigma_k(t),\quad x_k(t):= y_k \circ \sigma_k(t).
$$
Then, for every $k \in \cN,$ $u_k \in AC_L([a,b];U),$ where $L := {\rm Var}_{[0,1]} (\varphi),$ and $x_k:=x[\xb,u_k,v]$ is the Carath\'eodory solution of \eqref{E}-\eqref{IC} associated to $(u_k,v).$
Since the $y_k$ are equibounded, also the $x_k$ are equibounded.
Furthermore, for every $t\in [a,b],$
$$
|(u_k,x_k)(t)-(u,x)(t)| = |(\varphi,y_k)\circ \sigma_k(t) - (\varphi,y) \circ \sigma(t) | \to 0,
$$
since $(y_k)$ converges to  $y$ uniformly and $\sigma_k(t) \to \sigma(t).$
By the Dominated Convergence Theorem, we also have that $\|(u_k,x_k)- (u,x)\|_1 \to 0.$
Since the $u_k$ have equibounded variation, $x$ is a BV simple limit solution and so the ``only if'' part of Theorem \ref{gcpdgen1} is proved.
\epf


\subsection{Proof of the ``if'' part of Theorem \ref{gcpdgen1}}\label{Proofif}
Let us recall the statement of Theorem \ref{gcpdgen1}: {\em Let $\xb \in \cR^n$ and $(u,v)\in BV([a,b];U) \times L^1([a,b];V).$
 If $x:[a,b]\to\cR^n$ is a BV simple limit solution of \eqref{E}-\eqref{IC} corresponding to $(\xb,u,v),$ then $x$  is a g.c. solution.}

By hypothesis there exists $L>0$ and  a sequence  $(u_k)\subset AC_L([a,b];U)$ such that
\bel{appr1}
\lim(u_k,x_k)(t) = (u,x)(t),\quad \forall t\in [a,b],
\eeq
where, for every $k\in\cN,$ $x_k:=x[\xb,u_k,v]$ is the Carath\'eodory solution of \eqref{E}-\eqref{IC} corresponding to $\xb,u_k,v.$
For every $k\in\cN$ let us define a map $\sigma_k:[a,b]\to[0,1]$ by setting
\bel{hClock}
\sigma_k(t)= \frac{t-a+{\rm Var}_{[a,t]}(u_k)}{b-a+ {\rm Var}_{[a,b]}(u_k)}.
\eeq
Since $u_k$ is absolutely continuous,  $\sigma_k$ is absolutely continuous as well. Moreover, it is strictly increasing.   Let $\ubf_{0,k}:[0,1]\to [a,b]$ be the inverse of $\sigma_k$ and set $\ubf_k:=u_k\circ\ubf_{0,k}.$ The space-time controls
$(\ubf_{0,k}, \ubf_k)$ turn out to be equi-Lipschitz, with Lipschitz constants  bounded by $ b-a+L.$ Hence, passing to a subsequence (still denoted by  $(\ubf_{0_k}, \ubf_k)$) if necessary, there exists a space-time control $(\ubf_0,\ubf)$ such that
\bel{unif1}
(\ubf_{0,k}, \ubf_k)\to (\ubf_0,\ubf),\quad \text{uniformly on $[0,1].$}
\eeq
Let us set $y_k:=y[\ubf_{0_k}, \ubf_k;v],$ for all $k\in\cN,$  and $y:=y[\ubf_{0}, \ubf; v].$ In particular, one has
$x_k=y_k\circ\sigma_k,$ for all $k\in \cN.$


 We shall prove that $(\varphi_0,\varphi)$ is a g.c. of $u.$ Note first  that $\varphi_0$ is nondecreasing and $(L+b-a)$-Lipschitz continuous. Furthermore,  $\varphi(0)=a,$ $\varphi(1)=b.$
  Let us choose $t\in [a,b].$ Since $(\sigma_k(t)) \subset [0,1]$  for all $k\in\cN,$ there exists  a subsequence $(\sigma_{k^t_j} )_{j\in\cN}$ and $\sigma(t) \in [0,1]$ such that
\bel{sigmahj}
\lim \sigma_{k^t_j}(t)  = \sigma(t) .
\eeq
Hence
$$
\big(\varphi_0(\sigma(t)),\varphi(\sigma(t)) \big) = \lim (\varphi_{0,k^t_j}, \varphi_{k^t_j})\circ \sigma_{k^t_j}(t)=\lim (t,u_{k_j^t}(t)) = (t,u(t)).
$$
Since for every $t\in [a,b],$ we have found an $s\,(=\sigma(t))\in [0,1]$ such that
$
(t,u(t)) = (\varphi_0(s),\varphi(s)),
$
 it turns out that $(\varphi_0,\varphi)$   is  a  graph completion of $u $ and, moreover, the map $t\mapsto\sigma(t)$ is a clock corresponding to  $u$ and  $(\varphi_0,\varphi).$  In particular,
 \benl
 \sigma(t)\in\ubf_0^{\leftarrow}(t) .
 \eenl


Consider the set-valued g.c. solution $\hat x$ corresponding to $(\ubf_0,\ubf),$ namely
$$
\hat x(t):=x[\ubf_0,\ubf;v]=y\circ\ubf_0^{\leftarrow}(t),\quad\forall t\in [a,b],
$$
and define  the set valued map
$$
\hat u(t) := \ubf\circ\ubf_0^{\leftarrow}(t),\quad\forall t\in [a,b].
$$
So that, in particular, the map
\benl
t\mapsto (\tilde u,\tilde x)(t):=(\ubf,y) \circ\sigma(t) \in (\hat u(t),\hat x(t)), \qquad \forall t\in [a,b]
\eenl
is a single-valued graph completion.
Observe that
 \benl
 \begin{split}
{\rm d}\Big((u,x)(t),(\tilde u,\tilde x)(t)\Big) &\leq \\
{\rm d}\Big((u,x)(t),(u_{k^t_j},x_{k^t_j})(t)\Big)+
{\rm d} \Big( (\varphi_{k^t_j},y_{k^t_j})&\circ \sigma_{k^t_j}(t),(\varphi,y) \circ \sigma(t) \Big)  \to 0,
\end{split}
 \eenl
 where ${\rm d}$ denotes the Euclidean distance.
 Indeed, the first term tends to zero by hypothesis \eqref{appr1}; while the second term tends to zero by \eqref{unif1} and \eqref{sigmahj}.  Therefore
 $(u,x)(t)=(\tilde u,\tilde x)(t)$ for all $t\in [a,b],$
 so the ``if'' part of Theorem \ref{gcpdgen1} is proved.

\epf

\subsection{Proof of Theorem \ref{Existence2}}
Let us recall the statement of Theorem \ref{Existence2}:

{\em Let us assume that $U$ has the Whitney property.
 Then, for every initial value $\xb \in \cR^n$ and control pair $(u,v) \in BV([a,b];U)\times L^1([a,b];V)$ there exists  an associated BV simple limit solution of \eqref{E}-\eqref{IC}.}

\begin{proof}
 Let $\T \subset[a,b]$ be the subset of instants at which  $u$ is discontinuous. Since it is at most countable, one can write $\T =\{t_i:i\in\cN\}.$ Let us construct a graph completion $(\ubf_0,\ubf)$ of $u$ as follows.
\begin{enumerate}
\item For each $i\in\cN$ consider two curves
$
\tilde\gamma_{-}^i:[0,1]\to U,\quad \tilde\gamma_{+}^i:[0,1]\to U
$
verifying
$$
\tilde\gamma_{-}^i(0) = u(t_i-),\quad
\tilde\gamma_{-}^i(1) = u(t_i)=\tilde\gamma_{+}^i (0),\quad
 \tilde\gamma_{+}^i (1) = u(t_i+),
    $$
    and
    $$
{\rm Var}_{[0,1]} ({\tilde\gamma}_{-}^i)\leq M |u(t_i)-u(t_i-)|,\quad
     {\rm Var}_{[0,1]} ({\tilde\gamma}_{+}^i)\leq M |u(t_i+)-u(t_i)|,
$$
with $M\geq 1.$
Such curves exist because $U$ has the Whitney property.

\item
Consider the function $\sigma:[a,b]\to [0,1]$ defined by
$$
\sigma(t) := \frac{t-a+{\rm Var}_{[a,t]}(u)}{b-a+{\rm Var}_{[a,b]}(u)},
$$
which is well-defined given that $u$ is of bounded variation. Observe that $u$ is continuous, [left continuous, right-continuous] at a $t\in[a,b]$ if and only if  $\sigma$  is  continuous, [left continuous, right-continuous]  at $t.$ Set
$$
c^i_-:=\sigma(t_i-),\quad d_-^i=c_+^i:=\sigma(t_i),\quad d^i_+:=\sigma(t_i+).
$$
So $c^i_-\leq d^i_+=c_+^i\leq d^i_+ .$
Let $\hat \varphi_0$ be the unique increasing continuous function verifying $\hat\varphi_0 \circ \sigma(t)=t,$ for all $t\in [a,b].$

On $[0,1]$ define the function $\hat\varphi$ by setting
$$
\hat \varphi(s):=
\left\{
\ba{cl}
\gamma^i_-\left( \frac{s}{d^i_- -c^i_-} \right), &\quad \text{if, for some} \,i,\, c^i_- < d^i_- \text{and } s\in [c^i_-,d^i_-] ,\\
u(t_i-),&\quad \text{if, for some} \,i,\,\, s=c^i_-=d_-^i,\\

\gamma^i_-\left( \frac{s}{d^i_+ -c^i_+} \right), &\quad \text{if, for some} \,i,\, c^i_+ < d^i_+ \text{and } s\in [c^i_+,d^i_+] ,\\
u(t_i+),&\quad \text{if, for some} \,i,\,\, s=c^i_+=d_+^i,\\
\gamma^i_+\left( \frac{s}{d^i_+ -c^i_+} \right), &\quad  \text{if, for some}\,i,\,\, s\in [c^i_+,d^i_+],\ \text{and } c^i_+ < d^i_+,\\
\gamma^i_+(0),&\quad  \text{if, for some} \,i,\,\,s=c^i_+=d_+^i,\\
u\circ \hat\varphi_0(s),&\quad \text{if } s\in [0,1] \backslash \hat\varphi_0(\T)
\ea
\right.
$$
   \item The curve $(\hat \varphi_0,\hat \varphi):[0,1]\to [a,b]\times U$ is locally Lipschitz and has bounded variation: ${\rm Var}_{[0,1]}((\hat\varphi_0,\hat\varphi)) \leq 1+ (2M-1){\rm Var}_{[a,b]}(u).$ So, by first reparameterizing $(\hat \varphi_0,\hat \varphi)$ with arclenght $\xi$ and then rescaling the latter by   $s:=\xi\cdot {\rm Var}_{[0,1]}((\hat\varphi_0,\hat\varphi)),$ one obtains a Lipschitz continuous curve   $$ (\varphi_0,\varphi):[0,1]\to [a,b]\times U$$ defined on $[0,1]$ with  Lipschitz constant equal to ${\rm Var}_{[0,1]}((\hat\varphi_0,\hat\varphi)).$   In particular, $(\varphi_0,\varphi)$ is a graph completion of $u.$
   \end{enumerate}

Consider the map
$$
x:=y[\ubf_0,\ubf;v]\circ\sigma.
$$
The function $x$ is then  the  g.c. solution of \eqref{E}-\eqref{IC} corresponding to $(\ubf_0,\ubf),$ $v,$ and the clock $\sigma,$ namely $x=x[\ubf_0,\ubf;v]_\sigma.$ Therefore, by Theorem \ref{gcpdgen1}, $x$ is a BV simple limit solution of \eqref{E}-\eqref{IC}.
\end{proof}


\section{Likely developments}\label{conclSec}
We have mainly  analyzed the commutative case for $\L^1$ inputs and the noncommutative case for controls $u$ with bounded variation. In particular we have seen how some existing notions of solution (for discontinuous inputs $u$) can be embedded within the class of limit solutions or within some specific subclass.  Below  we briefly mention some few situations which we think might be worth studying through the notion of limit solution.

\subsection{Unbounded variation and noncommutative fields}  There exists  a definition  of solution   for the case where neither the commutativity hypothesis is assumed  nor the  controls $u$ have bounded variation (see e.g \cite{BreRam94}).  This solution  is based on a fibration  of the state space by means of the  leaves of the  ideal generated by the Lie brackets of order $\geq 2$ of the vector fields ${g_1,\dots,g_m}.$   Approximation results  of the flows generated by the brackets are used to prove that the so-called {\it looping controls} can be regarded as limit of ordinary controls. In particular, the {quotient} system resulting from a local factorization  turns out to be commutative. We conjecture  that the solutions resulting by this approach  can be proved to be (likely simple)  limit solutions.

\subsection{Continuous controls $u$}\label{continuousu} The theory of {\it rough paths} begun in 1994  by T.J. Lyons and nowadays variously  developed (see e.g.  \cite{Lyo94,Lyo02,Lej03}), is widely recognized as an effective construction to deal with {\it continuous}  inputs $u$ \footnote{ We are not aware of extension including an ordinary, measurable, control $v$ besides $\dot u$, maybe an issue not interesting the major fields of application of rough paths.}. It is impossible to give in the  restricted space of this subsection  even an approximative idea of the notion of solution to \eqref{E}  when the inputs $u$ are identified with  rough paths.  Let us only remind that the notion of rough paths was introduced as a {\it nonlinear} development of the concept of Young's integral  and successively has  become a powerful  tool in  the field of  stochastic differential equations. Let us also point out that two main notions are crucial within this theory:
\begin{itemize}
\item[1)] the  {\it $k$-iterated integrals}, namely
$$
I_{[s,t]}^{i_1,\dots,i_k} := \int_{s<t_1<\dots<t_k<t} {\rm d} u_{i_1}(t_1)\dots  {\rm d} u_{i_k}(t_k)
$$
which, in a sense, provide extra information besides the mere $1$-iterated integral $u(t)-u(s)$ to single out the proper solution\footnote{Notice, in particular, that $I_{[s,t]}^{i_1,i_2} -  I_{[s,t]}^{i_2,i_1}$ is twice the {\it area spanned} by the curve $(u_{i_2}, u_{i_2})$ during the interval $[t_1,t_2]$};
\item[2)] the pseudometric (on the  of inputs' and outputs' spaces) induced by the {\it $p$-variation}. Solutions are defined by means of  continuous extension of   input-output functionals with respect to these pseudo-metrics (the latter are rendered metrics by fixing initial points).  In particular, these solutions turn out to be unique and also  approachable  by suitable  Picard iterations.
\end{itemize}
Moreover, the use of rough paths involves Lie algebraic notions, incidentally demanding a sufficient amount  of regularity of vector fields to {\it  compensate  the roughness} of the paths.

A somehow different but connected approach,  concerning generalized solutions on the Heisenberg group  corresponding to H\"older-continuous inputs $u$, can be found in \cite{Mas08}.

Since the notions of solution given in the above-quoted literature were obtained as extensions of input-output maps in Banach  spaces,  they are likely  simple limit solutions in the sense proposed here. It would be perhaps instructive to study   sharper  connections between  these issues.

 \subsection{Bounded versus unbounded inputs} A  case to which the notion of limit solution could likely be extended  is that of systems where the $g_\alpha$ are $v$-dependent, namely systems of the form
\bel{Et}\dot x = {f}(t,x,u,v) +\sum_{\alpha=1}^m g_\alpha(x,u,v) \dot u_\alpha,\quad t\in [a,b].
\eeq
These kind of equations   are important in mechanical applications --e.g. when $u$ is a {\it shape parameter} and $v$ is a control representing an external force or torque- and in min-max control context where, for instance, the adjoint equations may contain a $v$-dependent term multiplied by an unbounded control, like in \eqref{Et} (see e.g. \cite{BarIsh89}).
Let us observe that the dependence of $g_\alpha$ on $v$ is much more critical than the $v$-dependence of $f,$ already in the case of $u$ with bounded variation, in that a simultaneous jump of $u$ and $v$ would make  the determination of the corresponding jump of $x$ quite delicate  (see e.g. \cite{Mil94,MotRam95,BetRam13}). The issue could  be actually extended  to more general equations of the form
$$
\dot x = \Phi(t,x,u,v,\dot u).
$$
See e.g. \cite{Mil96} for the latter general case, and \cite{BreRam10,RamSar00,PedTia09}  for the case when $\Phi$ is polynomial in $\dot u.$

\subsection{Nonsmooth vector fields}
Applications could require the weakening of the regularity assumptions on the vector fields $g_\alpha.$ For instance, if the latter were just locally Lipschitz continuous, conditions like {\it vanishing Lie brackets}  would still make sense in view of the results in \cite{RamSus07} or \cite{Sim96}, so the commutative case could be likely addressed. It should be noted, however, that the homeomorphism  utilized to represent limit solutions would be not differentiable, so further technical questions should be addressed to mimic the theory of the regular case.

\section*{Acknowledgements}
We wish to thank   P. Lamberti and  R. Monti (University of Padova) for a few bibliographical advice concerning Whitney sets,   G. De Marco (University of Padova) for helpful discussions on  functions which are pointwise limits of continuous maps, and M. Motta  (University of Padova) for several remarks on various issues of the paper. Finally, let us  acknowledge the anonymous referee for the   accurate reading of the manuscript and some quite stimulating inquires.

This work was partially supported by the European Union under the 7th Framework Programme FP7-PEOPLE-2010-ITN -  Grant agreement number 264735-SADCO, and the Fondazione CaRiPaRo Project (Italy)
``Nonlinear Partial Differential Equations: models, analysis, and
control-theoretic problems".


\end{document}